\documentclass[12pt, a4paper]{article}
\usepackage[T1]{fontenc}
\usepackage{amsmath,amsthm}
\usepackage{amssymb}
\usepackage{amscd}
\usepackage[latin1]{inputenc}
\usepackage{amsfonts}
\usepackage[matrix,arrow,curve]{xy}

\theoremstyle{plain}
\newtheorem{theo}{Theorem}[section]
\newtheorem{lem}[theo]{Lemma}
\newtheorem{prop}[theo]{Proposition}

\theoremstyle{definition}
\newtheorem{df}[theo]{Definition}
\newtheorem{theosd}[theo]{Theorem}
\newtheorem{rem}[theo]{Remark}

\newtheorem{propp}[theo]{Proposition}
\newtheorem{corro}[theo]{Corollary}

\title{$R$-equivalence on low degree complete intersections}
\author{Alena Pirutka}

\begin{document}

\maketitle
\date{}
\begin{abstract}
Let $k$ be a function field in one variable over $\mathbb C$ or  the field $\mathbb C((t))$. Let $X$ be a $k$-rationally simply connected variety defined over $k$. In this paper we show that  $R$-equivalence on rational points of $X$ is trivial and that the Chow group of zero-cycles of degree zero $A_0(X)$ is zero.  In particular, this  holds for a smooth complete intersection of $r$ hypersurfaces in $\mathbb P^n_k$ of respective degrees $d_1,\ldots,d_r$ with $\sum\limits_{i=1}^{r}d_i^2\leq n+1$.
\end{abstract}

\section{Introduction}
$\,\,\;$ Let $X$ be a projective variety over a field $k$. Two rational points $x_1, x_2$ of $X$ are called \textit{directly $R$-equivalent} if there is a morphism $f:\mathbb P^1_k\to X$ such that $x_1$ and $x_2$ belong to the image of $\mathbb P^1_k(k)$. This generates an equivalence relation called \textit{$R$-equivalence} \cite{Man}. The set of $R$-equivalence classes is denoted by $X(k)/R$. If $X(k)=\emptyset$ we set $X(k)/R=0$.

From the above definition, the study of $R$-equivalence on $X(k)$ is closely related to the study of rational curves on $X$, so that we need many rational curves on $X$. The following class of varieties sharing this property was introduced in the $1990$s by Koll\'ar, Miyaoka
and Mori, and independently by Campana.

\df{Let $k$ be a field of characteristic zero. A projective geometrically integral variety $X$ over $k$ is called \textit{rationally connected} if for any algebraically closed field $\Omega$ containing $k$,  two general $\Omega$-points $x_1, x_2$ of $X$ can be connected by a rational curve : there is a morphism $\mathbb P^1_{\Omega}\to X_{\Omega}$ such that its image contains $x_1$ and $x_2$.}

\rem{One can define rational connectedness by other  properties (\cite{Ko96}, Section IV.3).
For instance, if $k$ is uncountable, one may ask that the condition above is satisfied for \textit{any} two points $x_1, x_2\in X$.\\}

By a result of Campana \cite{Ca} and Koll\'ar-Miyaoka-Mori \cite{KMM}, smooth Fano varieties are rationally connected. In particular, a smooth complete intersection of $r$ hypersurfaces in $\mathbb P^n_k$ of respective degrees $d_1,\ldots,d_r$ with $\sum\limits_{i=1}^{r}d_i\leq n$ is rationally connected. Another important result about rationally connected varieties has been established by Graber, Harris and Starr \cite{GHS}. Let $k$ be  a function field in one variable over $\mathbb C$, that is, $k$ is the function field of a complex curve.  Graber, Harris and Starr prove that any smooth rationally connected variety over $k$ has a rational point.

One can see rationally connected varieties as an analogue of path connected spaces in topology. From this point of view, de Jong and Starr introduce the notion of rationally simply connected varieties as an algebro-geometric analogue of simply connected spaces. In this paper we use the following definition (see section $2$ for more details) :

\df\label{drsc}{Let $k$ be a field of characteristic zero. A projective geometrically integral variety $X$ over $k$ is called \textit{$k$-rationally simply connected} if for any sufficiently large integer $e$ there exists a geometrically irreducible component $M_{e,2}\subset \bar M_{0,2}(X,e)$ intersecting the open locus of irreducible curves $M_{0,2}(X,e)$ and such that the restriction of the evaluation morphism $$ev_2: M_{e,2}\to X\times X$$ is dominant with rationally connected general fiber.  \\}

Note that a $k$-rationally simply connected variety $X$ over a field $k$ is rationally connected as $X\times X$ is dominated by $M_{e,2}\cap M_{0,2}(X,e)$ from the definition above. This implies that two general points of $X$ over any algebraically closed field $\Omega\supset k$  can be connected by a rational curve.

By a recent result of de Jong and Starr \cite{dJS}, a smooth complete intersection $X$ of $r$ hypersurfaces in $\mathbb P^n_k$ of \,respective degrees $\,d_1,\ldots,d_r$ and of dimension at least $3$  is $k$-rationally simply connected if $\sum\limits_{i=1}^{r}d_i^2\leq n+1$.

Let us now recall the definition of $A_0(X)$. Denote $Z_0(X)$ the free abelian group generated by the closed points of $X$. \textit{The Chow group of degree $0$} is the quotient of the group $Z_0(X)$  by the subgroup generated by the $\pi_*(\mathrm {div}_C(g))$, where $\pi:C\to X$ is a proper morphism from a normal integral curve $C$, $g$ is a rational function on $C$ and $\mathrm{div}_C(g)$ is its divisor. It is denoted by $CH_0(X)$. If $X$ is projective, the degree map $Z_0(X)\to \mathbb Z$ which sends a closed point $x\in X$ to its degree $[k(x):k]$ induces a map $\mathrm{deg}: CH_0(X)\to \mathbb Z$ and we denote $A_0(X)$ its kernel.

At least in characteristic zero,
the set $X(k)/R$ and the group $A_{0}(X)$
are $k$-birational invariants of smooth projective $k$-varieties,
and they are reduced to
one element if $X$ is a projective space. Thus
  $X(k)/R=1$ and $A_{0}(X)=0$  if  $X$ is a smooth projective $k$-rational variety.

For $k$ a function field in one variable over $\mathbb C$ and for
$k=\mathbb
C((t))$,
one wonders whether a similar statement holds for arbitrary
smooth rationally connected $k$-varieties
(\cite{CT},  10.11 and 11.3).
This is established for some special classes of varieties (\cite{CT83}, \cite{CTSa}, \cite{CTSk}, \cite{CTSaSD}, \cite{Ma}). Most of these results hold as soon as $cd(k)\leq 1$ or even as soon as $k$ is a $C_1$-field.

In this paper, we prove the following result :

\theosd\label{main}{\textit{Let $k$ be either a function field in one variable over $\mathbb C$ or  the field $\mathbb C((t))$. Let $X$ be a $k$-rationally simply connected variety over $k$. Then
\begin{itemize}
\item[(i)] $X(k)/R=1$;
\item[(ii)] $A_0(X)=0$.
\end{itemize} }}

Combined with the theorem of de Jong and Starr, this gives :

\corro{\textit{ Let $k$ be either a function field in one variable over $\mathbb C$ or  the field $\mathbb C((t))$. Let $X$ be a smooth complete intersection of $r$ hypersurfaces in $\mathbb P^n_k$ of \,respective degrees $\,d_1,\ldots,d_r$.  Assume that $\sum\limits_{i=1}^{r}d_i^2\leq n+1.$ Then
\begin{itemize}
\item[(i)] $X(k)/R=1$;
\item[(ii)] $A_0(X)=0$.
\end{itemize} }}

The methods we use in the proof of the theorem apply more generally over a field $k$ of characteristic zero such that any rationally connected variety over $k$ has a rational point. As for the corollary, one can prove it in a simpler way for any $C_1$ field in the case $\sum d_i^2\leq n$ (see section $4$).\\

Note that one knows better results for smooth cubics and smooth intersections of two quadrics.  Let $k$ be either a function field in one variable over $\mathbb C$ or  the field $\mathbb C((t))$.

In the case of smooth cubic  hypersurfaces in $\mathbb P_k^n$ we have $X(k)/R=1$ if $n\geq 5$ (\cite{Ma}, 1.4). It follows that $A_{0}(X)=0$ if $n\geq 5$. In fact, we have $A_{0}(X)=0$ if $n\geq 3$. One can prove it by reduction to the case $n=3$. The latter case follows from the result on  geometrically rational $k$-surfaces (\cite{CT83}, Thm.A), obtained by $K$-theoretic methods.

In the case of smooth intersections of two quadrics in $\mathbb P_k^n$ we have  $X(k)/R=1$ if $n\ge 5$ (\cite{CTSaSD}, 3.27). Hence $A_{0}(X)=0$
 if $n \geq 5$. In fact, we also have  $A_{0}(X)=0$
 if $n=4$ as a particular case of \cite{CT83}, Thm.A.

It was also known that under the assumption of the theorem there is a bound $N=N(d_1,\ldots,d_r)$ such that $A_0(X)=0$ if $n\geq N$ (\cite{Pa}, 5.4). The bound $N$ here is defined recursively : $N(d_1,\ldots,d_r)=f(d_1,\ldots, d_r; N(d_1-1,\ldots,d_r-1))$ where the function $f$ grows rapidly with the degrees. For example, one can deduce that $N(3)=12$ and $N(4)=3+3^{12}$.\\

Theorem \ref{main} is inspired by the work of de Jong and Starr \cite{dJS} and we use their ideas in the proof. In section \ref{s2} we recall some notions about the moduli space of curves used in \cite{dJS} and we analyse the case when  we can deduce some information about $R$-equivalence on $X$ from the existence of a rational point on the moduli space. Next, in section \ref{s3}  we deduce Theorem \ref{main}.  In section \ref{s4} we  give an application to complete intersections.\\

\textbf{Acknowledgement :} I am very grateful to my advisor, Jean-Louis Colliot-Th\'el\`ene, for his suggestion to use the result of \cite{dJS}, for many useful discussions and the time that he generously gave me. I would like to express my gratitude to Jason Starr for his interest, for pointing out Proposition \ref{corc1}  and for allowing me to put his arguments into this paper. I also want to thank to Jean-Claude Douai and David Harari for very helpful discussions concerning the arguments in \ref{gerbes}.

\section{Rational points on a moduli space of curves}\label{s2}

\subsection{The moduli space $\bar M_{0,n}(X, d)$}
$\,\,\;$Let $X$ be a projective variety over a field $k$ of characteristic zero with an ample divisor $H$. Let $\bar k$ be an algebraic closure of $k$. While studying $R$-equivalence on  rational points of $X$ we need to work with a space parametrizing rational curves on $X$. We  fix the degree  of the curves we consider in order to have a space of finite type.

The space of rational curves of fixed degree on $X$  is not compact in general. One way to compactify it, due to Kontsevich, is to use stable curves.

\df{A \textit{stable} curve  over $X$ of degree $d$  with $n$ marked points is a datum $(C, p_1,\ldots,p_n, f)$ of
\begin{itemize}
\item [(i)] a proper geometrically connected reduced $k$-curve $C$ with only nodal singularities,
\item [(ii)] an ordered collection $p_1,\ldots, p_n$ of distinct smooth $k$-rational points of $C$,
\item [(iii)] a $k$-morphism $f:C\to X$ with $\mathrm{deg}_C \,f^*H=d$,
\end{itemize}
such that the stability condition is satisfied :
\begin{itemize}
\item [(iv)]  C has only finitely many $\bar k$-automorphisms fixing the points $p_1,\ldots, p_n$ and commuting with $f$.
\end{itemize}

We say that two stable curves $(C, p_1,\ldots,p_n, f)$ and $(C', p'_1,\ldots,p'_n, f')$ are \textit{isomorphic} if there exists an isomorphism $\phi: C\to C'$ such that $\phi(p_i)=p'_i$, $i=1,\ldots,n$ and $f'\circ \phi=f$.\\
}

We use the construction of Araujo and Koll\'ar \cite{AK} to parametrize stable curves.  They show that there exists a \textit{coarse moduli space} $\bar M_{0,n}(X,d)$ for all genus zero stable curves over $X$ of degree $d$ with $n$ marked points, which is a projective $k$-scheme (\cite{AK}, Thm. 50). Over $\mathbb C$ the construction was first given in \cite{FP}. In this paper, the authors fix a \textit{curve class} $\beta\in H^{2\mathrm{dim} X-2}(X,\mathbb Z)$ rather than the degree. They consider stable curves $f : C\to X$ such that the class of $f_*[C]$ is $\beta$. They prove that  there exists a coarse moduli space $\bar M_{0,n}(X,\beta )$ parametrizing all genus zero stable curves of class $\beta$ with $n$ marked points. The result in \cite{AK} holds over arbitrary, not necessarily algebraically closed field and, more generally, over a noetherian base.

We denote by $M_{0,n}(X,d)$ the open locus corresponding to irreducible curves and by $$ev_n: \bar M_{0,n}(X,d)\to \underbrace{X\times\ldots\times X}_n$$ the evaluation morphism which sends a stable curve to the image of its marked points.

When one says that $\bar M_{0,n}(X,d)$ is a coarse moduli space, it means that the following two conditions are satisfied :

\begin{itemize}
\item[(i)] there is a bijection of sets :
$$
\Phi:\left\{\begin{matrix} \mbox{isomorphism classes of}\\
\mbox{ genus zero stable curves over }\bar k\;\\f:C\to X_{\bar k}
\mbox{ with }n\mbox{ marked points, }\\
\mathrm{deg}_C \,f^*H=d\\
\end{matrix}\right\}\stackrel{\sim}{\to} \bar M_{0,n}(X,d)(\bar k);
$$
\item[(ii)] if  $\mathcal C\to S$ is a family of genus zero stable curves of degree $d$ with $n$ marked points, parametrized by a $k$-scheme $S$, then there exists a unique morphism $M_S: S\to \bar M_{0,n}(X,d)$ such that for every $s\in S(\bar k)$ we have $$M_S(s)=\Phi(\mathcal C_s).$$
\end{itemize}

In general, over nonclosed fields we do not have a bijection between isomorphism classes of stable curves and rational points of the corresponding moduli space, see \cite{AK} p.31. In particular, a $k$-point of $\bar M_{0,n}(X,d)$ does not  in general correspond  to a stable curve defined over $k$. So we need to use some more arguments if we want to study  $R$-equivalence on rational points of $X$ as we need to have some curves defined over $k$.

Let $P$ and $Q$ be two $k$-points of $X$. Suppose there exists  a stable curve $f:C\to X_{\bar k}$ over $\bar k$ with two marked points mapping to $P$ and $Q$, such that the corresponding point $\Phi(f)$ is a $k$-point of $\bar M_{0,2}(X,d)$. Even if we are not able to prove that $f$ is defined over $k$, using combinatorial arguments we will show   that $P$ and $Q$ are $R$-equivalent over $k$. Let us state the main result of this section.

\propp\label{g1}{\textit{Let $X$ be a projective variety over a field $k$ of characteristic zero.  Let $P$ and $Q$ be $k$-points of $X$. Let  $f:C\to X_{\bar k}$ be a stable curve over $\bar k$  of genus zero with two marked points mapping to $P$ and $Q$. Let $H$ be a fixed ample divisor on $X$ and let $d=\mathrm{deg}_C \,f^*H$.   If the corresponding point $\Phi(f)\in \bar M_{0,2}(X,d)$ is a $k$-point of $\bar M_{0,2}(X,d)$, then the points $P$ and $Q$ are $R$-equivalent over $k$. \\}}

\subsection{A lemma about graphs}

$\,\,\;$In the proof of Proposition \ref{g1}, starting from a curve defined over some extension of $k$ we will construct some curve defined over the base field $k$. In order to do this we will analyse the conjugates of the given curve. Let us first fix some notation.\\

Let $k$ be a field of characteristic zero. Let us fix an algebraic closure $\bar k$ of $k$.  Let $L\stackrel{i}{\hookrightarrow} \bar k$ be a finite Galois extension of $k$, and let $G=\mathrm{Aut}_k(L)$. For any $\sigma\in G$ we denote $\sigma^*: \mathrm{Spec}\, L\to \mathrm{Spec}\, L$ the induced morphism.
If $Y$ is an $L$-variety,  denote $\,^{\sigma}Y$ the base change   of $Y$ by $\sigma^*$ and  $\,^{\sigma}Y_{\bar k}$ the base change by $(i\circ\sigma)^*$. We denote the projection $\,^{\sigma}Y\to Y$ by  $\sigma^*$ too. If $f:Z\to Y$  is an $L$-morphism of $L$-varieties, then we denote $\,^{\sigma}f:\,^{\sigma}Z\to \,^{\sigma}Y$ and $\,^{\sigma}f_{\bar k}: \,^{\sigma}Z_{\bar k}\to \,^{\sigma}Y_{\bar k}$ the induced morphisms.

Note that if $Y\subset \mathbb P^n_L$ is a projective variety, then $\,^{\sigma}Y$ can
be obtained by applying $\sigma$ to each coefficient in the equations defining $Y$. Thus, if $Y$ is defined over $k$,  then the subvarieties $Y, \,^{\sigma}Y$ of  $\mathbb P^n_L$ are given by the same embedding for all $\sigma\in G$.
In this case the collection of morphisms $\{\sigma^*:Y\to Y\}_{\sigma\in G}$ defines a right action of $G$ on $Y$.  By Galois descent (\cite{BLR}, 6.2), if a subvariety $Z\subset Y$ is stable under this action of $G$, then $Z$ also is defined over $k$.\\

For lack of a suitable reference, let us next give a proof of the following lemma :

\lem\label{g2}{Let $C$ be a projective geometrically connected  curve  of arithmetic genus $p_a(C)=h^1(C,O_C)=0$ over a perfect field $k$. Assume $C$ has only nodal singularities. Then any two smooth $k$-points $a, b$ of $C$  are $R$-equivalent.}

\proof{For any field extension  $F$ of $k$ let us call an \textit{$F$-path joining $a$ and $b$} a closed $F$-subcurve $C'\subset C_F$ such that
\begin{itemize}
 \item [(i)] $C'=C'_1\cup C'_2\cup \ldots\cup C'_r$ where $C'_i$, $i=1,\ldots r$ are smooth $F$-rational curves;
\item[(ii)] $a\in C'_1$, $b\in C'_r$;
\item[(iii)]  if $1\leq i\leq r-1$, the intersection $C'_i\cap C'_{i+1}$ is an $F$-point and  the curves $C'_i$ and $C'_{i+j}$ do not intersect for $j>1$.
\end{itemize}

Since the  arithmetic genus of $C$ is zero, its geometric components are smooth rational curves over $\bar k$ intersecting transversally. As $C$ is geometrically connected, there exists a $\bar k$-path $C'$ joining $a$ and $b$. We may assume that $C'$ is an $L$-path for some finite Galois extension $L$ of $k$. Moreover, such a path is unique : if there were two different paths we would have a cycle formed by components of $C_{L}$, which is impossible as $p_a(C)=0$. We would like to find a $k$-path, thus we will achieve the proof.

Let us write $C'=C'_1\cup C'_2\cup \ldots\cup C'_r$, $a\in C'_1$, $b\in C'_r$. We will show that $C'$ comes from a $k$-path by base extension.   Let us take $\sigma\in \mathrm{Aut}_k(L)$. Then $\,^{\sigma}C'_1,\ldots, \,^{\sigma}C'_r$ is an $L$-path joining $a$ and $b$. Since such a path is unique, for every $i=1,\ldots r$ the components  $C'_i$ and $\,^{\sigma}C'_i$ of $C_{L}$  are equal and $C'_i\cap C'_{i+1}=\,^{\sigma}C'_i\cap \,^{\sigma}C'_{i+1}=\sigma(C'_i\cap C'_{i+1})$, $i=1,\ldots r$. This means that every component of the path $C'_i\subset C_L$  is stable over the action of  $\mathrm{Aut}_k(L)$ on $C_{L}$, hence it is defined over $k$, that is $C'_i=D_i\times_k L$,  for some $k$-curve $D_i\subset C$. By the same argument,  the intersection points of $D_i$ and $D_{i+1}$, $i=1,\ldots r-1$ are $k$-points. We deduce that $\{D_1,\ldots D_r\}$ is a $k$-path joining $a$ and $b$, so the points $a$ and $b$ are $R$-equivalent over $k$.\qed\\}

\subsection{Rational points on $\bar M_{0,2}(X,d)$}

$\;\,\;\;$ Let us now give the  proof of Proposition \ref{g1}.

We call $a$ and $b$ the marked points of $C$. We may assume that  $C$, $f$, $a$ and $b$ are defined over a finite Galois extension $L\stackrel{i}{\hookrightarrow} \bar k$ of $k$. That is, we may assume that $C$ is an $L$-curve, $a,b\in C(L)$ and that we have an $L$-morphism $f:C\to X_{L}$. Let us denote $T=\mathrm{Spec}\,L$. We view $L$ as a $k$-scheme and $f:C\to X\times_k T$ as a family of stable curves parametrized by $T$. Thus we have a moduli map $M_T: T=\mathrm{Spec}\,L\to \bar M_{0,2}(X,d)$ defined over $k$ and such that  for every $t\in T(\bar k)$ we have $$M_T(t)=\Phi(C_t)$$ where $f_t:C_t\to X_{\bar k}$ is the fibre of $f:C\to X\times_k T$ over $t$.

Note that $T\times_k\bar k=\prod\limits_{G} \mathrm{Spec}\,\bar k$ where the product is indexed by $G=\mathrm{Aut}_k(L)$
and the morphism $\prod\limits_{G} \bar k\to T=\mathrm{Spec}\,L$ is given by $(i\circ\sigma)^*$  on the corresponding component.
This implies that a $\bar k$-point $t\in T(\bar k)$ corresponds to some $\sigma\in G$ and  the morphism $f_t$ is the base
change by $(i\circ\sigma)^{*}$. Hence the morphism $f_t$ is the morphism $\,^{\sigma} f_{\bar k}:\,^{\sigma}C_{\bar k}\to X_{\bar k}$
and the marked points of $^{\sigma}C_{\bar k}$ are $\sigma(a)$ and $\sigma(b)$.

Since the curve $f_{\bar k}:C_{\bar k}\to X_{\bar k}$ corresponds to a $k$-point of $\bar M_{0,2}(X,d)$, we can factor $M_T$ as $$T=\mathrm{Spec}\, L\to \mathrm{Spec}\, k\stackrel{\Phi(f_{\bar k})}{\to} \bar M_{0,2}(X,d).$$
We thus see  that for every $t\in T(\bar k)$ the point $M_T(t)$ is the same point $\Phi(f_{\bar k})$ of $\bar M_{0,2}(X,d)$. 
 Hence for every $\sigma\in G$ the  curves $^{\sigma}C_{\bar k}$ and $C_{\bar k}$ are isomorphic as stable curves. This means that there exists a $\bar k$-morphism $\phi_{\sigma}:C_{\bar k}\to \,^{\sigma}C_{\bar k}$, such that $$\phi_{\sigma}(a)=\sigma(a), \;\phi_{\sigma}(b)=\sigma(b)\mbox{ and }^{\sigma}f_{\bar k}\circ\phi_{\sigma}=f_{\bar k}.$$

As a consequence, the proposition results from the following lemma.\\

\lem\label{g3}{Let $X$ be a projective variety over a perfect field $k$. Let $L$ be a finite Galois extension of $k$. Denote $G=\mathrm{Aut}_k(L)$. Let $P$ and $Q$ be $k$-points of $X$. Suppose we can find  an $L$-stable curve of genus zero $f:C\to X_{L}$  with two marked points $a,b\in C(L)$, satisfying the following conditions :
\begin{itemize}
\item[(i)] $f(a)=P$, $f(b)=Q$;
\item[(ii)] for every $\sigma\in G$ there exists a $\bar k$-morphism $ \phi_{\sigma}:C_{\bar k}\to \,^{\sigma}C_{\bar k}$ such that
$$\phi_{\sigma}(a)=\sigma(a), \;\phi_{\sigma}(b)=\sigma(b)
\mbox{ and }\,^{\sigma}f_{\bar k}\circ\phi_{\sigma}=f_{\bar k}.$$
\end{itemize}
Then the points $P$ and $Q$ are $R$-equivalent over $k$.}

\proof{ By lemma \ref{g2}, we have a unique $L$-path $\{C_1,\ldots ,C_m\}$ joining $a\in C_1(L)$ and $b\in C_m(L)$, where $C_i$ are irreducible components of $C$. We will use the curves $f(C_1),\ldots f(C_m)$ to show that $P$ and $Q$ are $R$-equivalent over $k$. Let us first show that these curves are defined over $k$ and not only over $L$.

For every $\sigma\in G$ we have an $L$-path $ \{\,^{\sigma}C_1,\ldots  ,\,^{\sigma}C_m\}$ joining $\sigma(a)\in \,^{\sigma}C_1(L)$ and $\sigma(b)\in \,^{\sigma}C_m(L)$. On the other hand,  $\{\phi_{\sigma}(C_{1, \bar k} ),\ldots ,\phi_{\sigma}(C_{m,\bar k})\}$ is a $\bar k$-path joining $\sigma(a)$ and $\sigma(b)$. Since the arithmetic genus of $^{\sigma}C_{\bar k}$ is zero, such a path is unique. That is, it coincides with the path $ \{\,^{\sigma}C_{1,\bar k},\ldots  ,\,^{\sigma}C_{m,\bar k}\}$. So we have $$\phi_{\sigma}(C_{i, \bar k})=\,^{\sigma}C_{i,\bar k},\;i=1,\ldots ,m.$$
Let us fix $1\leq i\leq m$. Denote the image $f(C_i)$ of $C_i$ in $X_L$ by $Z_i$. As $\,^{\sigma} f$ is a base change by $\sigma^*$, we have  the following commutative diagram :
$$\begin{CD}  \,^{\sigma}C_i @>\,^{\sigma}f>> \,^{\sigma}X\\
@VVV @VVV  @. \\
C_i @>f>> X
\end{CD}$$
We thus see that $\,^{\sigma}f(\,^{\sigma}C_i)=\,^{\sigma}Z_i$. Using base change by $i:L\to \bar k$ in the first line of the diagram above we obtain that $\,^{\sigma}f_{\bar k}(\,^{\sigma}C_{i,\bar k})=\,^{\sigma}Z_{i,\bar k}$.  On the other hand, since $\phi_{\sigma}(C_{i, \bar k})=\,^{\sigma}C_{i,\bar k}$ and $\,^{\sigma}f_{\bar k}\circ\phi_{\sigma}=f_{\bar k}$,  we have $\,^{\sigma}Z_{i,\bar k}=\,^{\sigma}f_{\bar k}(\,^{\sigma}C_{i,\bar k})= \,^{\sigma}f_{\bar k}(\phi_{\sigma}(C_{i, \bar k}))= f_{\bar k}(C_{i, \bar k})=Z_{i,\bar k}$. Since $\,^{\sigma}Z_{i}$ and $Z_{i}$ are $L$-subvarieties of $X_{L}$, we deduce that $\,^{\sigma}Z_{i}=Z_{i}$ for all $\sigma\in G$. By Galois descent, this means that the curve $Z_{i}$ is defined over $k$, that is, there exists a $k$-curve $D_i\subset X$ such that $Z_{i}=D_i\times_k L$.

In order to conclude the proof, we will show that the curve $D_i$ is a $k$-rational curve on $X$, that is, it is the image of some morphism from $\mathbb P^1_k$ to $X$, and that the point $f(C_i\cap C_{i+1})$ is the image of a $k$-point. Note that it may be not so obvious if $f(C_i\cap C_{i+1})$ is a singular point of $D_i$.

Let  $\tilde D_i\to D_i$ be the normalisation morphism. It induces an isomorphism over the smooth locus $D_i^{sm}$. Since $C_i$ is smooth, the morphism $f|_{C_i}:C_i\to D_i\times_k L$ extends to a morphism $f_i: C_i\to \tilde D_i\times_k L$ :
 $$\xymatrix{
 & \tilde D_i\times_k L\ar[d]&\\
 C_i\ar[r]^{f|_{C_i}} \ar[ur]^{f_i} & D_i\times_k L.
}
$$
This implies that $\tilde D_i$ is an $L$-rational curve. We have $\,^{\sigma}f_{i, \bar k}\circ\phi_{\sigma}=f_{i, \bar k},$ as this is true over a Zariski open subset $D_i^{sm}$. Moreover, for every $1\leq i\leq m-1$ we have $\phi_{\sigma}(C_i\cap C_{i+1})=\,^{\sigma}C_i\cap\,^{\sigma} C_{i+1}$. Using the same argument as above, we deduce that the point $f_i(C_i\cap C_{i+1})$ is a $k$-point of $\tilde D_i$.  This implies that $\tilde D_i$ is  a $k$-rational curve as it is $L$-rational and has a $k$-point. Moreover,  the point $f(C_i\cap C_{i+1})$ is a $k$-point of X as the image of $f_i(C_i\cap C_{i+1})$. Hence $P$ is $R$-equivalent to $f(C_1\cap C_2)$ as there is a rational curve $\tilde D_1\to X$ connecting them. By the same argument,  $f(C_{i-1}\cap C_{i})$ is $R$-equivalent to $f(C_{i}\cap C_{i+1})$ for all $1< i< m-1$ and $f(C_{m-1}\cap C_m)$ is $R$-equivalent to $Q$. Therefore $P$ and $Q$ are $R$-equivalent.    \qed\\}

\rem\label{gerbes}{If the cohomological dimension of the field $k$ is at most $1$, one can use more general arguments to prove Proposition \ref{g1}. In fact, one can consider the following stack $\mathcal G$ over the \'etale site Spec$\,k$ :
\begin{itemize} \item[(i)]the objects of $\mathcal G$ over the extension $L$ of $k$ are stable curves $D\to X_L$ over $L$  which become isomorphic over $\bar k$ to the
$\bar k$-curve $C\to X_{\bar k}$ corresponding to the point $\Phi(f)$ as in Proposition \ref{g1}.
 \item [(ii)] the morphisms between two stable curves over $L$ are $L$-isomorphisms.
\end{itemize}
One can also view  $\mathcal G$ as the fibre of the morphism $\overline{\mathcal M}_{0,2}(X,d)\to \bar M_{0,2}(X,d)$ over the point $\Phi(f)$, where $\overline{\mathcal M}_{0,2}(X,d)$ is the stack of all  genus zero stable curves over $X$ of degree $d$ with two marked points.

By this description, the objects of $\mathcal G$ exist locally (that is, over some finite extension of $k$). Moreover, any two such objects are locally isomorphic (that is, after taking some finite extension of $k$). A stack satisfying these two properties is called \textit{gerbe} (cf.\cite{Gi}). What we want to prove is the existence of objects over $k$, that is, that the curve $C\to X_{\bar k}$ is defined over $k$. A gerbe which has objects over $k$ is called \textit{neutral}.

By considering the automorphism groups  of objects of $\mathcal G$ one can associate a \textit{band} $\mathcal{L}$ to the gerbe $\mathcal G$, see \cite{Gi}, Ch.IV.  Next, one can define a cohomological set $H^2(\mathrm{Spec}k_{\text{\'et}}, \mathcal{L})$
 parametrizing all classes of gerbes whose associated band is $\mathcal L$. Note that in our case the automorphism group of objects of $\mathcal G$ is locally (that is, starting from some extension of $k$) a finite constant group consisting of the automorphisms of $C\to X_{\bar k}$ over $\bar k$. Under this assumption and also under the assumption that cd$\,k\leq 1$ one knows by \cite{DDE}, cor. 1.3, that all classes of gerbes in $H^2(\mathrm{Spec}k_{\text{\'et}}, \mathcal{L})$ are neutral. In particular, $\mathcal G$ is neutral.

From these arguments we conclude that there is a genus zero $k$-stable curve $f':C'\to X$ with two marked points $a, b\in C'(k)$ such that the images of $a$ and $b$ are respectively the points $P$ and $Q$. By lemma \ref{g2} we deduce that $a$ and $b$ are $R$-equivalent in $C'$, thus their images $P$ and $Q$ in $X$ are also $R$-equivalent.
}

\section{Proof of the theorem}\label{s3}
$\,\,\;$ In this section we use the previous arguments  to prove Theorem \ref{main}. Let $k$ be a function field in one variable over $\mathbb C$ or  the field $\mathbb C((t))$.  Let  $X$ be a $k$-rationally simply connected variety. In particular, $X$ is rationally connected.  Note that by the theorem of Graber, Harris and Starr \cite{GHS}, a smooth rationally connected variety over  a function field in one variable over $\mathbb C$ has a rational point, and the same result is  also known  over $k=\mathbb C((t))$ (cf. \cite{CT} 7.5). As any smooth projective variety equipped with a birational morphism to $X$ is still rationally connected, it has a rational point. This implies that $X(k)\neq\emptyset$.

 Let us fix a  sufficiently large integer $e$ and an irreducible component $M_{e,\,2}\subset \bar M_{0,2}(X,e)$ such that the restriction of the evaluation morphism $ev_2:M_{e,\, 2}\to X\times X$   is dominant with rationally connected general fibre.

Let $P$ and $Q$ be two $k$-points of $X$. A general strategy is the following.   We would like to apply \cite{GHS} and to deduce that there is a rational point in a fibre over  $(P,Q)$.  Then, by Proposition \ref{g1}, we deduce that $P$ and $Q$ are $R$-equivalent. But we only know that a general fibre of $ev_2$ is rationally connected. If $k=\mathbb C((t))$, this is sufficient as  $R$-equivalence classes are Zariski dense in this case by \cite{Ko99}. If $k$ is a function field in one variable over $\mathbb C$, our strategy will also work thanks to the result of Hogadi and Xu \cite{HX}. They prove that if we take a dominant proper morphism of $k$-varieties $h:Y\to Z$ such that $Z$ is smooth and the generic fibre of $h$ is rationally connected, then for every point $z\in Z$ there exists a subvariety of the fibre $Y_z$, defined over $k(z)$, which is geometrically irreducible and rationally connected.

Hence we can find a rationally connected $k$-subvariety $V$ in the fibre $ev_2^{-1}(P, Q)$. By what was recalled above, there is a rational point in $V$, hence in $ev_2^{-1}(P, Q)$. By the proposition \ref{g1}, the points $P$ and $Q$ are $R$-equivalent.  So we obtain $X(k)/R=1$ as $X(k)\neq\emptyset$.

Let us now prove that the group $A_0(X)$ is trivial. Pick  $x_0\in X(k)$. It is sufficient to prove that for every closed point $x\in X$ of degree $d$ we have that $x-dx_0$ is zero in $CH_0(X)$. Let us take a rational point $x'\in X_{k(x)}$ over a point $x$.
By the first part of the theorem, applied to $X_{k(x)}$, $x'$ is $R$-equivalent to $x_0$ over $k(x)$. Hence $x'-x_0$ is zero in $CH_0(X_{k(x)})$. Applying the push-forward by the morphism $p:X_{k(x)}\to X$,
we deduce that $x-dx_0$ is zero in $CH_0(X)$. This finishes the proof. \\ \qed

Note that one can prove the theorem without using the result of \cite{HX}. In fact, there is another way to prove the existence of a rational point in every fiber of the morphism $ev_2$, see \cite{Sta} p.25 for example. We thank M.Lieblich for this reference.  Let us sketch the argument here. More precisely, we will show the following :
\lem{ Let $k=\mathbb C(C)$ be  the function field of a (smooth) complex curve $C$. Let $Z$ and $T$ be projective $k$-varieties, with $T$ smooth. Let $f:Z\to T$ be a morphism with rationally connected general fibre. Then for \textit{every} $t\in T(k)$ there exists a rational point in the fibre $Z_t$.}
\proof{One can choose  proper models $ \mathcal T\to C$ and $F: \mathcal Z\to \mathcal T$ of $T$ and $ Z$ respectively with $\mathcal T$ smooth. We know that  any fibre of $F$ over some open set $U\subset\mathcal T$ is rationally connected.

The point $t\in T(k)$ corresponds to a section $s: C\to \mathcal T$. What we want is to find a section $C\to \mathcal Z\times_{\mathcal T} C$. One can view the image $s(C)$ in $\mathcal T$ as a component of a complete intersection $C'$ of hyperplane sections of $\mathcal T$ for some projective embedding. In fact, it is sufficient to take $\mathrm{dim}\,\mathcal T-1$ functions in the ideal of $s(C)$ in $\mathcal T$ generating this ideal over some open subset of $s(C)$. Moreover, one may assume that $C'$ is a special fibre  of a family $\mathcal C$ of hyperplane sections with general fibre a smooth curve intersecting $U$. After localization, we may also assume that $\mathcal C$ is  parametrized by $\mathbb C[[t]]$. Let $A$ be any affine open subset in $\mathcal C$ containing the generic point $\xi$ of $s(C)$.  We have the following  diagram :
$$\xymatrix{
 & & \mathcal Z\times_{\mathcal T}\mathrm{Spec}\,A\ar[d]^{F_A}\ar[r]&\mathcal Z\ar[d]^{F} &\\
 \xi\ar[r] & \mathrm{Spec}\,A\otimes_{\mathbb C[[t]]} \mathbb C\ar[r]\ar[d] & \mathrm{Spec}\,A\ar[d]\ar[r] &\mathcal T&\\
 & \mathrm \mathrm{Spec}\,\mathbb C \ar[r] & \mathrm{Spec}\,\mathbb C[[t]]. &
}
$$
Let $K=\mathbb C((t))$ and let $\bar K$ be an algebraic closure of $K$. By construction, the generic fibre of $F_{\bar K}: \mathcal Z\times_{\mathcal T} \bar K\to \mathrm{Spec}\,A \otimes_{\mathbb C[[t]]} \bar K $ is rationally connected.  By \cite{GHS}  we obtain a rational section of $F_{\bar K}$. As $\bar K$ is the union of the extensions $\mathbb C((t^{1/N}))$ for $N\in \mathbb N$,  we have a rational section for the morphism $\mathcal Z\times_{\mathcal T}\mathbb C[[t^{1/N}]]\to \mathrm{Spec}\,A \otimes_{\mathbb C[[t]]}  \mathbb C[[t^{1/N}]]$ for some $N$. By properness, this section extends to all codimension $1$ points of $\mathrm{Spec}\,A \otimes_{\mathbb C[[t]]}\mathbb C[[t^{1/N}]]$, in particular, to the point $\xi$ on the special fiber. This extends again  to give a section $C\to \mathcal Z\times_{\mathcal T} C$  as desired.\qed\\ }

\section{Proof of the corollary}\label{s4}

The following result is essentially contained in \cite{dJS}. We include the proof here as we need the precise statement over a field which is not algebraically closed.

\prop\label{cc}{Let $k$ be a field of characteristic zero. Let $X$ be a  smooth complete intersection of $r$ hypersurfaces in $\mathbb P^n_{k}$ of respective degrees $d_1,\ldots,d_r$ with $\sum\limits_{i=1}^{r}d_i^2\leq n+1$.  Suppose that $\mathrm{dim}\,X\geq 3$. Then for every $e\geq 2$  there exists a geometrically irreducible $k$-component $M_{e,2}\subset \bar M_{0,2}(X,e)$ such the restriction of the evaluation morphism $$ev_2:M_{e, 2}\to X\times X$$  is dominant with rationally connected generic fibre.}

\proof{ Let us first recall the construction of \cite{dJS} in the case $k=\mathbb C$. In this paper, the authors work with the space $\bar M_{0,2}(X,\beta)$ of \cite{FP} which parametrizes  stable curves of genus zero over $X$ of class $\beta\in H^{2\mathrm{dim\,X}-2}(X,\mathbb Z)$  with two marked points. Hovewer, as $\mathrm{dim}\,X\geq 3$, we know that  $H^{2\mathrm{dim\,X}-2}(X,\mathbb Z)=\mathbb Z\alpha$ where the degree of $\alpha$ equals to $1$ (\cite{V}, 13.25).  Thus we can replace $\beta$ by its degree $e$ and work with the space  $\bar M_{0,2}(X, e)$ as in \cite{AK}.

In \cite{dJS}, de Jong and Starr prove that for every integer $e\geq 2$ there exists an irreducible component $M_{e,2}\subset \bar M_{0,2}(X,e)$ such that the restriction of the evaluation morphism $ev_2:M_{e, 2}\to X\times X$ is dominant with rationally connected generic fibre. We will specify more precisely how they get the component $M_{e, 2}$. It will follow from their construction  that $M_{e,2}$ is in fact the unique component satisfying the above property.
The construction of $M_{e,2}$ is the following :
\begin{enumerate}
\item One first shows that there exists a \emph{unique} irreducible component $M_{1,1}\subset \bar M_{0,1}(X, 1)$ such that the restriction of the evaluation $ev_1|_{M_{1,1}}: M_{1,1}\to X$ is dominant (\cite{dJS}, 1.7).
\item The component $M_{1,\,0}\subset \bar M_{0,\,0}(X,1)$ is constructed as the image of $M_{1, 1}$ under the morphism $\bar M_{0,1}(X, 1)\to \bar M_{0,\,0}(X, 1)$  forgetting the marked point. Then one constructs the component of higher degree $M_{e,0}$ as the \emph{unique} component of $\bar M_{0,\,0}(X, e)$ which intersects the subvariety of $\bar M_{0,\,0}(X, e)$ parametrizing a degree $e$ cover of the smooth, free curve parametrized by $M_{1,\,0}$ (\cite{dJS}, 3.3).
\item The component $M_{e,2}\subset \bar M_{0,2}(X,e)$ is the \emph{unique} component such that its image under the morphism $\bar M_{0,2}(X, e)\to \bar M_{0,\,0}(X, e)$, which forgets about the marked points, is $M_{e, 0}$.
\end{enumerate}

Let us now consider the general case. Let $\bar k$ be an algebraic closure of $k$. As $k$ is of finite type over $\mathbb Q$, we may assume that $\bar k\subset\mathbb C$.  Since the decomposition into geometrically irreducible components does not depend on which algebraically closed field we choose, by the first step above  there exists a unique irreducible component $M_{1,1}\subset \bar M_{0,1}(X_{\bar k}, 1)$ such that the restriction of the evaluation $ev_1|_{M_{1,1}}$ is dominant. As this component is unique, it is  defined over $k$. Hence, from the construction above,  the component $M_{e,2}$ is also defined over $k$, which completes the proof. \qed\\}

Let us now prove the corollary. Let $X$ be a  smooth complete intersection of $r$ hypersurfaces in $\mathbb P^n_{k}$ of respective degrees $d_1,\ldots,d_r$ with $\sum\limits_{i=1}^{r}d_i^2\leq n+1$. Thus, if dim$\,X=1$ then $X$ is a line and the corollary is obvious. If dim$\,X=2$ then  $X$ is a quadric surface in $\mathbb P^3_k$. We have that $X$ is birational to $\mathbb P^2_k$ as it has a $k$-point and the corollary follows.   If dim$\,X\geq 3$, we have that $X$  is $k$-rationally simply connected  by the Proposition \ref{cc}. If $k$ is a function field in one variable over $\mathbb C$ or  the field $\mathbb C((t))$ we have  $X(k)/R=1$ and $A_0(X)=0$ by the theorem \ref{main}. This finishes the proof of the corollary.\\

\rem{The next argument, due to Jason Starr, gives a simpler way to prove the corollary in the case $\sum d_i^2\leq n$. More precisely, we have~:
\prop\label{corc1}{Let $k$ be a $C_1$ field. Let $X\stackrel{i}{\hookrightarrow} \mathbb P_k^n$ be the vanishing set of $r$ polynomials $f_1,\ldots f_r$ of respective degrees $d_1,\ldots d_r$. If $\sum d_i^2\leq n$ then any two points $x_1, x_2\in X(k)$ can be joined by two lines defined over $k$ : there is a point $x\in X(k)$ such that $l(x,x_i)\subset X$, $i=1,2$, where $l(x,x_i)$ denote the line through $x$ and $x_i$. }
\proof{We may assume that $x_1=(1:0:\ldots :0)$ and $x_2=(0:1:0:\ldots :0)$ via the embedding $i$. The question is thus to find a point $x=(x_0:\ldots:x_n)$ with coordinates in $k$ such that $$\left\{
                                                           \begin{array}{ll}
                                                             f_i(tx_0+s,tx_1,\ldots tx_n)=0\\
                                                            f_i(tx_0, tx_1+s,\ldots tx_n)=0,
                                                           \end{array}
                                                         \right.
 i=1,\ldots r.$$ As $x_1, x_2$ are in $X(k)$ these conditions are satisfied for $t=0$. Thus we may assume $t=1$. Writing $f_i(x_0+s,x_1,\ldots x_n)=\sum\limits_{j=0}^{d_i}P_j^i(x_0,\ldots x_n)s^j$ with deg$P_j^i=d_i-j$ we see that each equation $f_i(x_0+s,x_1,\ldots x_n)=0$ gives us $d_i$ conditions on $x_0,\ldots x_n$ of degrees $1, \ldots d_i$. By the same argument, each equation $f_i(x_0, x_1+s,\ldots x_n)=0$ gives $d_i-1$ conditions of degrees $1, \ldots d_i-1$ as we know from the previous equation that we have no term of degree zero. The sum of the degrees of all these conditions on $x_0,\ldots x_n$ is $\sum\limits_{i=0}^{r} d_i^2$. As $\sum\limits_{i=0}^{r} d_i^2\leq n$ by Tsen-Lang theorem we can find a solution over $k$, which finishes the proof.  \qed\\}}

\end{document}